\documentclass[leqno]{article}
\usepackage{amsmath, amsthm}

\theoremstyle{plain}
\newtheorem{theorem}{Theorem}[section]

\theoremstyle{definition}
\newtheorem{definition}[theorem]{Definition}

\numberwithin{equation}{section}
\date{}

\begin{document}
\title{\textbf{Some Characterizations of Slant Helices in the Euclidean Space} $\mathbf E^n$}

\author{\text{Ahmad T. Ali}\\
Mathematics Department\\
 Faculty of Science, Al-Azhar University\\
 Nasr City, 11448, Cairo, Egypt\\
E-mail: \textit{atali71@yahoo.com}\\
\vspace*{1cm}\\
\text{Melih Turgut}\footnote{Corresponding author.} \\
Department of Mathematics, \\
Buca Educational Faculty, Dokuz Eyl\"{u}l University,\\
35160 Buca, Izmir, Turkey\\
E-mail: \textit{melih.turgut@gmail.com, melih.turgut@ogr.deu.edu.tr}}

\maketitle
\begin{abstract}
In this work, notion of a slant helix is extended to space E$^n$. Necessary and sufficient conditions to be a slant helix in the Euclidean $n-$space are presented. Moreover, we express some integral characterizations of such curves in terms of curvature functions.\\

\textbf{M.S.C. 2000}: 53A04
\end{abstract}

\section{Introduction and Statement of Results}
Inclined curves or so-called general helices are well-known curves in the classical differential geometry of space curves \cite{mp} and we refer to the reader for recent works on this type of curves \cite{gl1,sc}. Recently, Izumiya and Takeuchi have introduced the concept of slant helix in Euclidean 3-space E$^3$ saying that the normal lines makes a constant angle with a fixed direction \cite{izu}. They characterize a slant helix if and only if the function
\begin{equation}\label{eqma}
\dfrac{\kappa^2}{(\kappa^2+\tau^2)^{3/2}}\Big(\dfrac{\tau}{\kappa}\Big)^{\prime}
\end{equation}
is constant. In the same space, spherical images, the tangent and the binormal indicatrix and some characterizations of such curves are presented by \cite{ly}. With the notion of a slant helix, similar works are treated by the researchers, see \cite{ali, al2, ey, mmhk, ms}.

In this work, we consider the generalization of the concept of a slant helix in the Euclidean n-space E$^n$.\\

Let $\alpha:I\subset R \rightarrow E^n$ be an arbitrary curve in E$^n$. Recall that the curve $\alpha$ is said to be of unit speed (or parameterized by arc-length function $s$) if $\langle\alpha'(s),\alpha'(s)\rangle=1$, where $\langle,\rangle$ is the standard scalar product in Euclidean space E$^n$ given by
\begin{center}
$\langle X,Y\rangle=\sum_{i=1}^{n}\,x_i\,y_i,$
\end{center}
for each $X=(x_1,\ldots,x_n)$, $Y=(y_1,\ldots,y_n)\in E^n$.

Let $\{V_{1}(s),\ldots,V_{n}(s)\}$ be the moving frame along $\alpha$, where the vectors $V_{i}$ are mutually orthogonal vectors satisfying $\langle V_{i}, V_{i}\rangle=1$.
The Frenet equations for $\alpha$ are given by (\cite{gl1})
\[
\left[
\begin{array}{c}
V_{1}^{\prime } \\
V_{2}^{\prime } \\
V_{3}^{\prime } \\
\vdots  \\
V_{n-1}^{\prime } \\
V_{n}^{\prime }
\end{array}
\right] =\left[
\begin{array}{ccccccc}
0 & \kappa _{1} & 0 & 0 & \cdots  & 0 & 0 \\
-\kappa _{1} & 0 & \kappa _{2} & 0 & \cdots  & 0 & 0 \\
0 & -\kappa _{2} & 0 & \kappa _{3} & \cdots  & 0 & 0 \\
\vdots  & \vdots  & \vdots  & \vdots  & \ddots  & \vdots  & \vdots  \\
0 & 0 & 0 & 0 & \cdots  & 0 & \kappa _{n-1} \\
0 & 0 & 0 & 0 & \cdots  & -\kappa _{n-1} & 0
\end{array}
\right] \left[
\begin{array}{c}
V_{1} \\
V_{2} \\
V_{3} \\
V_{4} \\
V_{5} \\
V_{6}
\end{array}
\right].
\]

Recall the functions $\kappa_i(s)$ are called the i-th curvatures of $\alpha$. If $\kappa_{n-1}(s)=0$ for any $s\in I$, then $V_{n}(s)$ is a constant vector $V$ and the curve $\alpha$ lies in a $(n-1)$-dimensional affine subspace orthogonal to $V$, which is isometric to the Euclidean $(n-1)$-space E$^{n-1}$. We will assume throughout this work that all the  curvatures satisfy $\kappa_i(s)\not=0$ for any $s\in I$, $1\leq i\leq n-1$. Here, recall that a regular curve with constant Frenet curvatures is called a $W-$curve \cite{ps}.

\begin{definition} A unit speed curve $\alpha:I\rightarrow E^n$ is called slant helix if its unit principal normal $V_2$ makes a constant angle with a fixed direction $U$.
\end{definition}
Our main result in this work is the following characterization of slant helices in Euclidean $n$-space E$^n$.

\begin{theorem}\label{th-main} Let $\alpha:I\rightarrow E^n$ be a unit speed curve in E$^n$. Define the functions
\begin{equation}\label{u211}
\begin{array}{ll}
G_1=\int\kappa_{1}(s) ds,\ G_2=1,\ G_3=\dfrac{\kappa_1}{\kappa_2}G_1,\
G_{i}=\dfrac{1}{\kappa_{i-1}}\Big[\kappa_{i-2}G_{i-2}+G_{i-1}^{\prime}\Big],
\end{array}
\end{equation}
where $4\leq i\leq n$.
Then $\alpha$ is a slant helix if and only if the function
\begin{equation}\label{u21}
\sum_{i=1}^{n}\,G_i^2=C
\end{equation}
is constant and non-zero. Moreover, the constant $C=\sec^2\theta$, being $\theta$ the angle that makes $V_2$ with the fixed direction $U$ that determines $\alpha$.
\end{theorem}

This theorem generalizes in arbitrary dimensions what happens for $n=3$, namely: if $n=3$, (\ref{u21}) writes
$$1+\Big(1+\dfrac{\kappa_1^2}{\kappa_2^2}\Big)G_1^2=C.$$
It is easy to prove that: this equation is equivalent to
$$\dfrac{\kappa_1^2}{(\kappa_1^2+\kappa_2^2)^{3/2}}\Big(\dfrac{\kappa_2}{\kappa_1}\Big)^{\prime}=
\dfrac{1}{\sqrt{C-1}}$$
equation (\ref{eqma}) where $C\neq 1.$

\section{Proof of Theorem \ref{th-main}}
Let $\alpha$ be a unit speed curve in E$^n$. Assume that $\alpha$ is a slant helix curve. Let $U$ be the direction with which $V_2$ makes a constant angle $\theta$ and, without loss of generality, we suppose that $\langle U,U\rangle=1$.  Consider the differentiable functions $a_i$, $1\leq i\leq n$,
\begin{equation}\label{u3}
U=\sum_{i=1}^{n}\,a_i(s)\,V_i(s),\ \ s\in I,
\end{equation}
that is,
$$a_i=\langle V_i,U\rangle,\ 1\leq i\leq n.$$
Then the function  $a_2(s)=\langle V_2(s),U\rangle$ is constant, and it agrees with $\cos\theta$ as follows:
\begin{equation}\label{u311}
a_2(s)=\langle V_2,U\rangle=\cos\theta
\end{equation}
for any $s$. Because the vector field $U$ is constant, a differentiation in (\ref{u3}) together (\ref{u211}) gives the following system of ordinary differential equation:
\begin{equation}\label{u5}
\left\{\begin{array}{ll}
a_1^{\prime}-\kappa_1 a_2&=0\\
\kappa_1 a_1-\kappa_2 a_3&=0\\
a_{i}'+\kappa_{i-1} a_{i-1}-\kappa_{i} a_{i+1}&=0,\,\,\,\,\,3\leq i \leq n-1\\
a_n'+\kappa_{n-1} a_{n-1} &=0.
\end{array}\right.
\end{equation}
Let us define the functions $G_i=G_i(s)$ as follows
\begin{equation}\label{u51}
a_i(s)=G_i(s)\,a_2, \  1\leq i\leq n.
\end{equation}
We point out that $a_2\not=0$: on the contrary, (\ref{u51}) gives $a_i=0$, for $1\leq i\leq n$ and so, $U=0$, which is a contradiction. Since, the first $n$-equations in (\ref{u5}) lead to
\begin{equation}\label{u6}
\left\{
\begin{array}{ll}
&G_1=\int\kappa_{1}(s) ds\\
&G_2=1\\
&G_3=\dfrac{\kappa_1}{\kappa_2}G_1\\
&G_{i}=\dfrac{1}{\kappa_{i}}\Big[\kappa_{i-2}G_{i-2}+G_{i-1}^{\prime}\Big],\,\,\,\,\,4\leq i \leq n.\\
\end{array}\right.
\end{equation}
The last equation of (\ref{u5}) leads to the following condition;
\begin{equation}\label{u62}
G_n'+\kappa_{n-1} G_{n-1}=0.
\end{equation}
We do the change of variables:
$$
t(s)=\int^s\kappa_{n-1}(u) du,\hspace*{.5cm}\frac{dt}{ds}=\kappa_{n-1}(s).
$$
In particular, and from the last equation of (\ref{u6}), we have
$$
G_{n-1}'(t)=G_n(t)-\Big(\frac{\kappa_{n-2}(t)}{\kappa_{n-1}(t)}\Big)G_{n-2}(t).
$$
As a consequence, if $\alpha$ is a slant helix, substituting the equation (\ref{u62}) to the last equation, we express
$$G_n''(t)+G_n(t)=\frac{\kappa_{n-2}(t)G_{n-2}(t)}{\kappa_{n-1}(t)}.$$
By the method of variation of parameters, the general solution of this equation is obtained
\begin{equation}\label{u9}
G_n(t)=\Big(A-\int\frac{\kappa_{n-2}(t)G_{n-2}(t)}{\kappa_{n-1}(t)}\sin{t}\,dt\Big)\cos{t}+
\Big(B+\int\frac{\kappa_{n-2}(t)G_{n-2}(t)}{\kappa_{n-1}(t)}\cos{t}\,dt\Big)\sin{t},
\end{equation}
where $A$ and $B$ are arbitrary constants. Then (\ref{u9}) takes the following form
\begin{equation}\label{u10}
\begin{array}{ll}
G_n(s)=&\Big(A-\int\Big[\kappa_{n-2}(s)G_{n-2}(s)\sin{\int\kappa_{n-1}(s)ds}\Big]ds\Big)\cos{\int\kappa_{n-1}(s)ds}\\
&+
\Big(B+\int\Big[\kappa_{n-2}(s)G_{n-2}(s)\cos{\int\kappa_{n-1}(s)ds}\Big]ds\Big)\sin{\int\kappa_{n-1}(s)ds}.
\end{array}
\end{equation}
From (\ref{u62}), the function $G_{n-1}$ is given by
\begin{equation}\label{u11}
\begin{array}{ll}
G_{n-1}(s)=&\Big(A-\int\Big[\kappa_{n-2}(s)G_{n-2}(s)\sin{\int\kappa_{n-1}(s)ds}\Big]ds\Big)\sin{\int\kappa_{n-1}(s)ds}\\
&-\Big(B+\int\Big[\kappa_{n-2}(s)G_{n-2}(s)\cos{\int\kappa_{n-1}(s)ds}\Big]ds\Big)\cos{\int\kappa_{n-1}(s)ds}.
\end{array}
\end{equation}
From Equation (\ref{u6}), we have
\begin{eqnarray*}
\sum_{i=1}^{n-2}G_i G_i^{\prime}&=&G_1G_1'+G_2G_2'+\sum_{i=3}^{n-2}G_i G_i^{\prime}\\
&=&\kappa_1G_1+\sum_{i=3}^{n-2}G_i\Big[\kappa_iG_{i+1}-\kappa_{i-1}G_{i-1}\Big]\\
&=&\kappa_1G_1+\sum_{i=3}^{n-2}\Big[\kappa_iG_iG_{i+1}-\kappa_{i-1}G_{i-1}G_i\Big]\\
&=&\kappa_1G_1+\kappa_{n-2}G_{n-2}G_{n-1}-\kappa_{2}G_{2}G_3\\
&=&\kappa_{n-2}G_{n-2}G_{n-1}.
\end{eqnarray*}
Substituting (\ref{u11}) to the above equation and integrating it, we have:
\begin{equation}\label{u15}
\begin{array}{ll}
\sum_{i=1}^{n-2}G_{i}^{2}&=C-\Big(A-\int\Big[\kappa_{n-2}(s)G_{n-2}(s)\sin{\int\kappa_{n-1}ds}\Big]ds\Big)^2\\
&-
\Big(B+\int\Big[\kappa_{n-2}(s)G_{n-2}(s)\cos{\int\kappa_{n-1}ds}\Big]ds\Big)^2,
\end{array}
\end{equation}
where $C$ is a constant of integration. Using equations (\ref{u10}) and (\ref{u11}), we have
\begin{equation}\label{u151}
\begin{array}{ll}
G_n^2+G_{n-1}^2&=\Big(A-\int\Big[\kappa_{n-2}(s)G_{n-2}(s)\sin{\int\kappa_{n-1}ds}\Big]ds\Big)^2\\
&+
\Big(B+\int\Big[\kappa_{n-2}(s)G_{n-2}(s)\cos{\int\kappa_{n-1}ds}\Big]ds\Big)^2,
\end{array}
\end{equation}
It follows from  (\ref{u15}) and (\ref{u151}) that
$$\sum_{i=3}^{n}G_{i}^{2}=C.$$
Moreover this constant $C$ can be calculated as follows.
From (\ref{u51}), together the $(n-2)$-equations (\ref{u6}), we have
$$C=\sum_{i=1}^{n}\,G_i^2=\dfrac{1}{a_2^2}\sum_{i=1}^{n}\,a_i^2=\dfrac{1}{a_2^2}=\sec^2\theta,$$
where we have used (\ref{u21}) and the fact that $U$ is a unit vector field.

We do the converse of Theorem. Assume that the condition (\ref{u6}) is satisfied for a curve $\alpha$. Let $\theta\in R$ be so that $C=\sec^2\theta$. Define the unit vector $U$  by
$$U=\cos\theta\Big[\sum_{i=1}^{n}\,G_i\,V_i\Big].$$
By taking account (\ref{u6}), a differentiation of $U$ gives that $\dfrac{dU}{ds}=0$, which it means that $U$ is a constant vector field. On the other hand, the scalar product between the unit tangent vector field $V_2$ with $U$ is
$$
\langle V_2(s),U\rangle=\cos\theta.
$$
Thus, $\alpha$ is a slant helix in the space E$^{n}$.

As a direct consequence of the proof, we generalize theorem \ref{th-main} in Minkowski space for timelike curves and give an another theorem which characterizes slant helices with constant curvatures.

\begin{theorem} Let E$_1^n$ be the Minkowski n-dimensional space  and let $\alpha:I\rightarrow E_1^n$ be a unit speed timelike curve.
Then $\alpha$ is a slant helix if and only if
the function $\sum_{i=i}^{n}\,G_i^2$ is constant, where the functions $G_i$ are defined as in (\ref{u211}).
\end{theorem}
\begin{proof} The proof carries the same steps as above and we omit the details. We only point out that the fact that $\alpha$ is timelike means that $V_1(s)=\alpha'(s)$ is a timelike vector field. The other $V_i$ in the Frenet frame, $2\leq i\leq n$, are unit spacelike vectors and so, the second equation in Frenet equations changes to $V_2'=\kappa_1V_1+\kappa_2V_3$ (for details of Frenet equations see \cite{ehi}).
\end{proof}
\begin{theorem}
There are no slant helices with constant and non-zero curvatures ($W-$slant helices, i.e.) in the space E$^n$.
\end{theorem}
\begin{proof}
Let us suppose a slant helix with constant and non-zero curvatures. Then the equations in (\ref{u5}) and (\ref{u6}) hold. Since, we easily have for odd $i$, $G_{i}=\delta _{i}s$, where $\delta_{i}\in R$ and for even $i,$ $G_{i}=\delta _{i}$. Then, we form
\[
\sum\limits_{i=1}^{n}G_{i}^{2}=(\delta _{1}s)^{2}+\delta _{2}^{2}+(\delta _{3}s)^{2}+\delta _{4}^{2}+...
\]
and it is easy to say that $\sum\limits_{i=1}^{n}G_{i}^{2}$ is nowhere
constant. By the theorem \ref{th-main}, we arrive at that there does not exist a slant helix with constant and non-zero curvatures in the space E$^n$.
\end{proof}
\section{Further Characterizations of Slant Helices in E$^n$}
In this section we present new characterizations of slant helix in E$^n$. The first one is a consequence of Theorem \ref{th-main}.

\begin{theorem} \label{th-2} Let $\alpha:I\subset R\rightarrow E^n$ be a unit speed curve in Euclidean space E$^n$. Then $\alpha$ is a slant helix if and only if there exists a $C^2$-function $G_{n}(s)$ such that
\begin{equation}\label{u24}
G_{n}=\dfrac{1}{\kappa_{n-1}}\Big[\kappa_{n-2}G_{n-2}+G_{n-1}^{\prime}\Big],\,\,\,
\dfrac{dG_{n}}{ds}=-\kappa_{n-1}(s)G_{n-1}(s),
\end{equation}
where
$$G_1=\int\kappa_{1}(s) ds, G_2=1, G_3=\dfrac{\kappa_1}{\kappa_2}G_1, G_{i}=\dfrac{1}{\kappa_{i-1}}\Big[\kappa_{i-2}G_{i-2}+G_{i-1}^{\prime}\Big],\ 4\leq i\leq n-1.$$
\end{theorem}
\begin{proof}
Let now assume that $\alpha$ is a slant helix. By using Theorem \ref{th-main} and by differentiation the (constant) function  given in
(\ref{u21}), we obtain

\begin{eqnarray*}
0&=&\sum_{i=i}^{n}G_i\,G_i^{\prime}\\
&=&G_1\kappa_1+G_3\Big(\kappa_3G_4-\kappa_2G_2\Big)+G_4\Big(\kappa_4G_5-\kappa_3G_3\Big)+...\\
& &+G_{n-1}\Big(\kappa_{n-1}G_{n}-\kappa_{n-2}G_{n-2}\Big)+G_{n}G_{n}^{\prime}\\
&=&G_{n}\Big(G_{n}'+\kappa_{n-1}G_{n-1}\Big).
\end{eqnarray*}
This shows (\ref{u24}). Conversely, if (\ref{u24}) holds, we define a  vector field $U$ by
$$U=\cos\theta\Big[\sum_{i=1}^{n}\,G_i\,V_i\Big].$$
By the Frenet equations, $\dfrac{dU}{ds}=0$, and so, $U$ is constant. On the other hand,  $\langle V_2(s),U\rangle=\cos\theta$ is constant, and this means that $\alpha$ is a slant helix.
\end{proof}

We end  giving an integral characterization of a slant helix.

\begin{theorem} Let  $\alpha:I\subset R\rightarrow E^n$ be a unit speed curve in Euclidean space E$^n$. Then $\alpha$ is a slant helix if and only if the following condition is satisfied
\begin{equation}\label{u244}
\begin{array}{ll}
G_{n-1}(s)=&\Big(A-\int\Big[\kappa_{n-2}G_{n-2}\sin{\int\kappa_{n-1}ds}\Big]ds\Big)\sin{\int^s\kappa_{n-1}(u)du}\\
&-\Big(B+\int\Big[\kappa_{n-2}G_{n-2}\cos{\int\kappa_{n-1}ds}\Big]ds\Big)\cos{\int^s\kappa_{n-1}(u)du}.
\end{array}
\end{equation}
for some constants $A$ and $B$.
\end{theorem}
\begin{proof}
Suppose that $\alpha$ is a slant helix. By using  Theorem \ref{th-2}, let define $m(s)$ and $n(s)$ by
$$\phi(s)=\int^s\kappa_{n-1}(u)du,$$
\begin{equation}\label{u26}
\begin{array}{ll}
m(s)=&G_{n}(s)\cos\phi+G_{n-1}(s)\sin\phi+\int\,\kappa_{n-2}G_{n-2}\sin\phi\,ds,\\
n(s)=&G_{n}(s)\sin\phi-G_{n-1}(s)\cos\phi-\int\,\kappa_{n-2}G_{n-2}\cos\phi\,ds.
\end{array}
\end{equation}
If we differentiate equations (\ref{u26}) with respect to $s$ and taking into  account of (\ref{u244}) and (\ref{u24}), we
obtain  $\dfrac{dm}{ds}=0$ and $\dfrac{dn}{ds}=0$. Therefore, there exist constants $A$ and $B$ such that $m(s)=A$ and $n(s)=B$.
By substituting  into (\ref{u26}) and solving the resulting equations for $G_{n-1}(s)$, we get
$$
G_{n-1}(s)=\Big(A-\int\,\kappa_{n-2}G_{n-2}\sin{\phi}\,ds\Big)\sin{\phi}-
\Big(B+\int\,\kappa_{n-2}G_{n-2}\cos{\phi}\,ds\Big)\cos{\phi}.
$$

Conversely, suppose that (\ref{u244}) holds. In order to apply Theorem \ref{th-2}, we define $G_n(s)$ by
$$G_n(s)=\Big(A-\int\,\kappa_{n-2}G_{n-2}\sin{\phi}\,ds\Big)\cos{\phi}+
\Big(B+\int\,\kappa_{n-2}G_{n-2}\cos{\phi}\,ds\Big)\sin{\phi}.$$
with $\phi(s)=\int^s\kappa_{n-1}(u) du$. A direct differentiation of (\ref{u244}) gives
$$
G_{n-1}^{\prime}=\kappa_{n-1}G_{n}-\kappa_{n-2}G_{n-2}.
$$
This shows the left condition in (\ref{u24}). Moreover, a straightforward computation leads to
$G_{n}'(s)=-\kappa_{n-1}G_{n-1}$, which finishes the proof.
\end{proof}

ACKNOWLEDGEMENTS: The second author would like to thank T\"{u}bitak-Bideb for their financial supports during his Ph.D. studies.

\end{document}